\documentclass[11pt]{amsart}%%%%%,leqno
\usepackage{amsmath,amssymb, graphicx, amscd,latexsym,here}
\makeatletter
\newtheorem{Theorem}{Theorem}%[section]
\newtheorem{MainTheorem}[Theorem]{Main Theorem}%[section]
\newtheorem{Lemma}[Theorem]{Lemma}
\newtheorem{Corollary}[Theorem]{Corollary}
\newtheorem{Proposition}[Theorem]{Proposition}

\newtheorem{Remark}[Theorem]{Remark}

\newcommand{\eps}{\varepsilon}

\newcommand\la{\lambda}
\newcommand\vphi{\varphi}

\newcommand\al{\alpha}

\newcommand\ga{\gamma}

\newcommand\BC{ {\mathbb C}}
\newcommand\BN{ {\mathbb  N}}

\newcommand\BR{ {\mathbb  R}}
\newcommand\BP{ {\mathbb  P}}

\newcommand\bfs{\mbox {\bf  s}}

\newcommand\bfu{\mbox {\bf  u}}

\newcommand\bfv{\mbox {\bf  v}}

\newcommand\bfw{\mbox {\bf  w}}

\newcommand\bfz{\mbox {\bf  z}}

\newcommand\bfa{\mbox {\bf  a}}
\newcommand\bfb{\mbox {\bf  b}}

\newcommand\grad{\rm{grad}\/}

\newcommand\id{\rm{id}}

\newcommand\lcm{\rm{lcm}\/}

\newcommand\inv{^{-1}}

\def\mapright#1{\smash{\mathop{\longrightarrow}\limits^{{#1}}}}%\def\maprightt#1#2{\smash{\mathop{\longrightarrow}\limits^{#1}}}
\def\mapdown#1{\Big\downarrow\rlap{$\vcenter{\hbox{$#1$}}$}}

\def\inv{^{-1}}

\begin{document}
\title[On Mixed Brieskorn variety 
%of plane curves and non-reduced degeneration
%{\it Draft: \today}
]
{On Mixed Brieskorn variety 
}

\author
%Normally smooth divisors  {\it Draft: \today}
[M. Oka ]
{Mutsuo Oka}
\dedicatory{ Dedicated to Professor   A. Libgober for his 60th birthday}
% \end{center}
\address{\vtop{
\hbox{Department of Mathematics}
\hbox{Tokyo  University of Science}
\hbox{26 Wakamiya-cho, Shinjuku-ku}
\hbox{Tokyo 162-8601}
\hbox{\rm{E-mail}: {\rm oka@rs.kagu.tus.ac.jp}}
}}

%\today
%\thanks{ {Donated to Professor  Anatoly Libgober for his 60th birthday}
\begin{abstract}
Let $f_{\bfa,\bfb}(\bfz,\bar\bfz)=z_1^{a_1+b_1}\bar z_1^{b_1}+\cdots+z_n^{a_n+b_n}\bar
 z_n^{b_n}$  be a polar weighted  homogeneous mixed polynomial
 with $a_j>0,b_j\ge 0$, $j=1,\dots, n$
 and 
let  $f_{\bfa}(\bfz)=z_1^{a_1}+\cdots+z_n^{a_n}$ be the associated
 weighted homogeneous polynomial.
 Consider the corresponding link variety
 $K_{\bfa,\bfb}=f_{\bfa,\bfb}\inv(0)\cap S^{2n-1}$ and
 $K_{\bfa}=f_{\bfa}\inv(0)\cap S^{2n-1}$.
 Ruas-Seade-Verjovsky \cite{R-S-V} proved that
 the Milnor fibrations of $f_{\bfa,\bfb}$ and $f_{\bfa}$ are
 topologically equivalent and 
 the mixed link  $K_{\bfa,\bfb}$ is homeomorphic to
 the complex link $K_{\bfa}$. We will prove that they are
$C^\infty$ equivalent
 and two links are diffeomorphic. We show the same assertion for 
$ f(\bfz,\bar\bfz)=z_1^{a_1+b_1}\bar
 z_1^{b_1}z_2+\cdots+z_{n-1}^{a_{n-1}+b_{n-1}}\bar
 z_{n-1}^{b_{n-1}}z_n+z_n^{a_n+b_n}\bar z_n^{b_n}$ and its associated
 polynomial
 $ g(\bfz)=z_1^{a_1}z_2+\cdots+ z_{n-1}^{a_{n-1}}z_n+z_n^{a_n}$.
\end{abstract}
\maketitle

\maketitle

\section{Introduction}
%Let $f_{\bfa,\bfb}(\bfz,\bar\bfz)=z_1^{a_1+b_1}\bar z_1^{b_1}+\cdots+z_n^{a_n+b_n}\bar
% z_n^{b_n}$ and  $f_{\bfa}(\bfz)=z_1^{a_1}+\cdots+z_n^{a_n}$
% with $a_j>0,b_j\ge 0$ for each $j=1,\dots, n$.
We consider the mixed polar weighted homogeneous polynomial
$f_{\bfa,\bfb}(\bfz,\bar\bfz)$ and its associated polynomial $f_{\bfa}(\bfz)$:
\[\begin{split}
 &f_{\bfa,\bfb}(\bfz,\bar\bfz)=
z_1^{a_1+b_1}\bar z_1^{b_1}+\cdots+z_n^{a_n+b_n}\bar
   z_n^{b_n},
 \,\, f_{\bfa}(\bfz)=z_1^{a_1}+\cdots+z_n^{a_n}
   \end{split}
\] with $ a_j>0,\,b_j\ge 0,\,j=1,\dots,n$ and $\bfa=(a_1,\dots, a_n)$
and $\bfb=(b_1,\dots, b_n)$
and  we consider
the mixed Brieskorn variety
\[
 V_{\bfa,\bfb}:=\{\bfz\in \BC^n\,|\, f_{\bfa,\bfb}(\bfz,\bar
   \bfz)=0\},\, V_{\bfa}=\{\bfz\in \BC^n\,|\, f_{\bfa}(\bfz)=0\}.
\]
Put $K_{\bfa,\bfb,r}=V_{\bfa,\bfb}\cap S_r^{2n-1}$
and
$K_{\bfa,r}=V_{\bfa}\cap S_r^{2n-1}$.
Note that $f_{\bfa,\bfb}$ and $f_{\bfa}$ have the same
polar weights $P={}^t(p_1,\dots,p_n)$. Let $d=\lcm(a_1,\dots, a_n)$.
% (the least common multiple). 
Then $p_j=d/a_j$ for $j=1,\dots, n$.
Let us denote the associated $\BC^*$ action by
$t\circ z:=(z_1t^{p_1},\dots, z_nt^{p_n})$ for $t\in \BC^*$.
Then we have
\[\begin{split}
& f_{\bfa,\bfb}(\rho\circ \bfz,\overline{\rho\circ
 \bfz})=\rho^d\,f_{\bfa,\bfb}(\bfz,\bar\bfz),\,\,\rho\in S^1\subset \BC^*\\
& f_{\bfa}(t\circ \bfz)=t^{d}\,f_{\bfa}(\bfz),\,\,\qquad t\in \BC^*.
   \end{split}
\]
Note  that $f_{\bfa,\bfb}$ is also a radially weighted homogeneous
polynomial
and it defines a local and a  global Milnor fibrations which are
homotopically equivalent (\cite{OkaPolar,OkaMix}).

The Milnor fibration of
$f_{\bfa,\bfb}$:
\[
 f/|f|:\,S_{r}^{2n-1}\setminus\{K_{\bfa,\bfb,r}\}\to S^1
\]
does not depend on the radius $r$ and it
is topologically equivalent
to that of the complex polynomial $f_a(\bfz)=z_1^{a_1}+\cdots+z_n^{a_n}$
(\cite{R-S-V, OkaPolar}):
\[
  f:\BC^n\setminus V_{\bfa,\bfb}\to \BC^{*}.
 \]
Thus hereafter we put $K_{\bfa,\bfb}=K_{\bfa,\bfb,1}$
and $K_{\bfa}=K_{\bfa,r}$.
 Consider the homeomorphism: $\eta: \BC^n\to \BC^n$ defined by
 $\eta(\bfz)=(w_1,\dots, w_n)$ with
 $w_j=z_j|z_j|^{2b_j/a_j}$,\, $j=1,\dots, n$. Note that $\eta$ 
  preserves the values of $f_{\bfa,\bfb}$ and $f_{\bfa}$ and $\eta$ is 
 $S^1$-action equivariant. That is,
 \[
\eta(\rho\circ \bfz)=\rho\circ\eta(\bfz),\quad  f_{\bfa}(\eta(\bfa))=f_{\bfa,\bfb}(\bfz).
 \]
 Then  $\eta$ gives a homeomorphism of the two fibrations
$f_{\bfa,\bfb}:\BC^n\setminus V_{\bfa,\bfb}\to \BC^*$
and $f_{\bfa}:\BC^n\setminus V_{\bfa}\to \BC^*$ 
and a homeomorphism of the two hypersurfaces 
 $V_{\bfa,\bfb},\, V_{\bfa}$.
 Thus the following diagrams are commutative.
 \[
\begin{matrix}
 V_{\bfa,\bfb}&\subset& \BC^n&\supset& \BC^n\setminus V_{\bfa,\bfb}&\mapright{f_{\bfa,\bfb}}&\BC^*\\
 \mapdown{\eta}&&\mapdown{\eta}&&\mapdown{\eta}&&\mapdown{\id}\\
 V_{\bfa}&\subset& \BC^n&\supset& \BC^n\setminus V_{\bfa}&\mapright{f_{\bfa}}&\BC^*
 \end{matrix}
 \]
The homeomorphism
$\vphi: (S^{2n-1},K_{\bfa,\bfb})\to (S^{2n-1},K_{\bfa})$ is given
with a little modification of $\eta$
by
\[\begin{split}
 \vphi:\, & S^{2n-1}\to S^{2n-1},\,\vphi(\bfz)=\psi(\eta(\bfz))
   \end{split}
   \]
   where $\psi$ is the ``normalization'' mapping which is defined by
   $\bfw\mapsto r(\bfw)\circ\bfw$ where a positive real number $r(\bfw)$
   is defined by the equality: $\|r(\bfw)\circ \bfw\|=1$.
%and   $\bfw=\eta(\bfz)$.
It is easy to see that $\vphi$ gives also  a topological equivalence of the Milnor
fibrations:

\vspace{.2cm}
\[
 \begin{matrix}
  S^{2n-1}\setminus{K_{\bfa,\bfb}}&\mapright{f_{\bfa,\bfb}/|f_{\bfa,\bfb}|}&S^1\\
  \mapdown{\vphi}&&\mapdown{\id}\\
   S^{2n-1}\setminus{K_{\bfa}}&\mapright{f_{\bfa}/|f_{\bfa}|}&S^1
  \end{matrix}
  \] 
 Unfortunately we  observe  that neither $\eta$ nor $\vphi$ are 
differentiable
on the coordinate planes $z_j=0$.

The purpose of this note is to show that $\vphi$ (and $\eta$  also) %in fact
can  be replaced by a
diffeomorphism
$\vphi'$ which is isotopic to the identity map of the sphere $S^{2n-1}$
(Theorem \ref{MainTheorem}).
In \S 3, we prove the similar assertion for the polar weighted
homogeneous polynomial
\[
 f(\bfz,\bar\bfz)=z_1^{a_1+b_1}\bar
 z_1^{b_1}z_2+\cdots+z_{n-1}^{a_{n-1}+b_{n-1}}\bar
 z_{n-1}^{b_{n-1}}z_n+z_n^{a_n+b_n}\bar z_n^{b_n}\]
and its associated polynomial
$g(\bfz)=z_1^{a_1}z_2+\cdots+ z_{n-1}^{a_{n-1}}z_n+z_n^{a_n}$
(Main-Theorem-bis \ref{MainTheoremBis}).

Throughout this paper, we use the same notations as in \cite{OkaPolar,OkaMix}.

\section{Canonical family and the construction of an isotopy}
We consider the following
 mixed Brieskorn polynomial and its associated weighted
homogeneous polynomial in the sense of \cite{OkaPolar}:
\[\begin{split}
 &f_{\bfa,\bfb}(\bfz,\bar\bfz)=
z_1^{a_1+b_1}\bar z_1^{b_1}+\cdots+z_n^{a_n+b_n}\bar
   z_n^{b_n},
 \,\, f_{\bfa}(\bfz)=z_1^{a_1}+\cdots+z_n^{a_n}
   \end{split}
\]
First we consider the linear  family which connect two polynomials:
\[\begin{split}
 f_t(\bfz,\bar\bfz)&:=(1-t)f_{\bfa,\bfb}(\bfz,\bar\bfz)+t f_{\bfa}(\bfz)\\
&=(1-t)(z_1^{a_1+b_1}\bar z_1^{b_1}+\cdots+z_n^{a_n+b_n}\bar
   z_n^{b_n})+t(z_1^{a_1}+\cdots+z_n^{a_n})\\
   &=\sum_{j=1}^n \,z_j^{a_j}\left(t+(1-t)|z_j|^{2b_j}\right)
\end{split}
\]
for $0\le t\le 1$ and put $V_t=f_t\inv(0)$.
First we observe  that $f_0=f_{\bfa,\bfb}$ and $f_1=f_{\bfa}$ and 
$f_t,\,0\le t\le 1$ is a family of mixed polynomials which are polar weighted 
by the same weight
$P={}^t(p_1,\dots, p_n)$
where 
$p_j=\lcm(a_1,\dots, a_n)/a_j$, $j=1,\dots, n$,
though $f_t$ is not radially weighted homogeneous for $t\ne 0,1$.
Recall that the polar action is given as 
\[
 (\la,\bfz)\in S^1\times \BC^{*n},\,\, (\la,\bfz)\mapsto
 (z_1\la^{p_1},\dots,z_n\la^{p_n})\in \BC^{*n}.
\]
(In \cite{OkaPolar,OkaMix}, we have assumed a polar weighted homogeneous
polynomial
is also radially weighted homogeneous. In this paper, we do not  assume
the radial weighted homogeneity.)
%Put $V_t:=f_t\inv(0)
%\subset \BC^n$.
The first key assertion is:
\begin{Lemma}\label{key1}
The mixed polynomial $f_t(\bfz,\bar\bfz):\BC^n\to \BC$ has
%variety $V_t$ has 
a unique singularity at
 the origin and $V_t\setminus\{O\}$
and $f_t\inv(\eta)$ is mixed non-singular for any
 $t,\,0\le t\le 1$ and $\eta\ne 0$.
\end{Lemma}
\begin{proof}
Assume that %$\bfw\in V_{\tau}\setminus\{O\}$
 $\bfw\in \BC^n\setminus\{O\}$
is a mixed-singular  point of $f_\tau$
for some $0\le\tau\le 1$.
 As $V_0, V_1$ are mixed non-singular outside of the origin
 (\cite{OkaPolar}),
 we may
 assume that  $0<\tau<1$.
 We will show that this gives a contradiction.
By Proposition 1 of \cite{OkaPolar},
there exists a complex number $\la$ with $|\la|=1$ so that
\[
 \overline{df_{\tau}(\bfw,\bar\bfw)}=\la\,\bar d f_{\tau}(\bfw,\bar\bfw)
\]
where 
\[
 df_{\tau}=(\frac{\partial f_{\tau}}{\partial z_1},\dots,
\frac{\partial f_{\tau}}{\partial z_n}),\,
\bar d f_{\tau}=(\frac{\partial f_{\tau}}{\partial \bar z_1},\dots,
\frac{\partial f_{\tau}}{\partial \bar z_n}).
\]
This implies that
% for each
%$j,\, j=1,\dots, n$, we have the equality:
\[\begin{split}
 (a_j+b_j)\bar w_j^{a_j+b_j-1}w_j^{b_j}(1-{\tau})&+a_j\bar w_j^{a_{j-1}}\,{\tau}
=
b_j\,w_j^{a_j+b_j}\bar w_j^{b_j-1}(1-{\tau})\,\la,\\
%&\qquad\,\qquad\text{for}\,\,j=1,\dots, n.
\end{split}
\] for $j=1,\dots, n$.
Multiplying $\bar w_j$, we get the equality: 
\begin{eqnarray}
 \bar w_j^{a_j}\left\{
 (a_j+b_j)|w_j|^{2b_j}(1-{\tau})+a_j{\tau}\right\}=w_j^{a_j}|w_j|^{2b_j} b_j\la\,
 (1-{\tau}).
 %,\,j=1,\dots, n
\end{eqnarray} %for $j=1,\dots,n$.
Denote  the left side and the right side of (1)
 by $L(1)$ and $R(1)$ respectively. 
Then we have the inequality:
 \[\begin{split}
&  |L(1)|\ge |w_j|^{a_j+2b_j}(a_j+b_j)(1-{\tau})\ge|w_j|^{a_j+2b_j}b_j (1-{\tau})
    =|R(1)|.
   % &\qquad\qquad\,0<\tau<1
    \end{split}
  \]
 where the equality hold if and only if $w_j=0$.
As $L(1)=R(1)$ and $\tau\ne 1$, we must have $w_j=0$ for $j=1,\dots,n$.
That is, $\bfw=O$. This give a contradiction to the assumption
 $\bfw\ne O$.
\end{proof}
The second key observation is:
\begin{Lemma} \label{key2}
For any $t,\,0\le t\le 1$
 and for any $r>0$, the sphere $S_r^{2n-1}$ intersects $V_t$ transversely.
 \end{Lemma}
 \begin{proof}
  %This is the key assertion for the main result.
 As $f_t$ is not radially
  weighted homogeneous, the assertion is not obvious.
  Assume that the intersection is not transverse
  at $\bfw\in V_\tau\cap S_r^{2n-1}$ with $r=\|\bfw\|$.
 As $V_0,V_1$ are radially weighted homogeneous and 
any sphere $S_r^{2n-1}$ is transverse to them, we have that  $0<\tau<1$.
  As we have seen in the above  Lemma 1 that $V_\tau$ is mixed non-singular,
  the tangent space has the real codimension two.
 Let $f_t=g_t+ih_t$ where $g_t$ and $h_t$ be the real and the imiginary
  part of $f_t$, considering $g_t,h_t$ as functions of $2n$ variables
  $x_1,y_1,\dots, x_n,y_n$ with $z_j=x_j+iy_j$.
  The the tangent vectors are
  those vectors which are
  transverse (in the real Eucledean space $\BR^{2n}$) to
  the real gradient vectors $\grad\,g_\tau(\bfw)$ and $\grad\,h_\tau(\bfw)$
  at $\bfw$. 
 Non-transversality implies that three vectors
  $\bfw,\grad\,g_\tau(\bfw),\, \grad\,h_\tau(\bfw)$ are linearly dependent over $\BR$ at $\bfw$.
  As the latter two vectors are linearly independent,
  the tangent space $T_{\bfw}V_\tau$ is a subspace of
  the tangent space
  $T_{\bfw}S_r^{2n-1}$.
  We will show that this is impossible by showing the existence of a
  tangent vector $\bfu\in T_{\bfw}V_t$ which is not tangent to the
  sphere $S_r^{2n-1}.$
  We use the following simple assertion.
  \begin{Proposition}\label{prop}
Let $a_j>0,\,b_j\ge 0$ be fixed integers and let $\tau$ be a
   positive real number with $1> \tau>0$.
   For any  fixed $z\in \BC^*$ and $0\le\tau\le 1$, the function
   \[
    \psi_j(s,z):=|z|^{a_j}s^{a_j}\left(\tau+(1-\tau)|z|^{2b_j}s^{2b_j}\right)
   \]
   is a
   strictly monotone increasing function of $s$ on
the half line $\BR^+=\{s\,|\,s>0\}$ %for     any fixed  $\tau$
   and $j=1,\dots,n$.
   \end{Proposition}

\begin{proof}
  The proof is an easy calculus of the differential
$\frac{d\psi_j}{ds}$.
  In fact, the assertion follows from the strict positivity
  $\frac{d\psi_j}{ds}(s,z)>0$.
 \end{proof}

  Now we continue the proof of  Lemma \ref{key2}.
  Put $I=\{j\,|\, w_j\ne 0\}$.
  We may assume that $I=\{1,\dots,n\}$
  and $\bfw\in V_\tau\cap \BC^{*n}$ for simplicity.
  (Otherwise, we work in $\BC^{*I}$.)
  We use Proposition \ref{prop}  and the inverse function theorem
  to the function:
  \[
   r:\BR^+\to \BR^+,\, s\mapsto r(s)=\frac{\psi_j(s,w_j)}{|w_j|^{a_j}(\tau+(1-\tau)|w_j|^{2b_j})},\,r(1)=1
  \]
  to find   real-valued real analytic functions
  $\vphi_j(r,w_j)$ of the variable $r>0$ such that
  $\vphi_j(1,w_j)=1$ and 
  \[\begin{split}
   \psi_j(\vphi_j(r,w_j),w_j)&=r\,|w_j|^{a_j}(\tau+(1-\tau)|w_j|^{2b_j}),\,\,
    % &\qquad\qquad
     j=1,\dots, r.
     \end{split}
  \]
  This is equivalent to
  \begin{eqnarray}
   \vphi(r,w_j)^{a_j}&(\tau+(1-\tau)|w_j|^{2b_j}\vphi(r,w_j)^{2b_j})\\
   &=
  r\,|w_j|^{a_j}(\tau+(1-\tau)|w_j|^{2b_j})\notag.
  \end{eqnarray}
  By the monotone increasing property of the function $\psi_j(s,w_j)$,
  $\vphi_j(r,w_j)$ is also monotone increasing function of $r$ on 
the half line $0<r<\infty$.
   Now we consider the real analytic path
  $\xi(r)$ which is defined on $0<r<\infty$ by
  \[
  \xi(r)=(\eta_1(r,w_j),\dots, \eta_n(r,w_j)),\,\,
\text{where}\,\, \eta_j(r,w_j)=\vphi_j(r,w_j)w_j,\,\,j=1,\dots,n.
  \]
  As we have
  \[
   \arg\,\eta_j(r,w_j)^{a_j}(\tau+(1-\tau)|\eta_j(r)|^{2b_j})=
%\arg\,\eta_j(r)^{a_j}=
  \arg\,w_j^{a_j},
\]
  it is easy to observe that
  \[\begin{split}
   f_{\tau}(\xi(r),\bar\xi(r))&=\sum_{j=1}^n
     \eta_j(r)^{a_j}(\tau+(1-\tau)|\eta_j(r)|^{2b_j})\\
     &\overset{(2)} =\sum_{j=1}^n r\,w_j^{a_j}(\tau+(1-\tau|w_j|^{2b_j}))\\
     &=r\, f_\tau(\bfw,\bar\bfw)\equiv 0.
     \end{split}
  \]
  Thus $r\mapsto\xi(r),\,0\le r\le \infty$ is a curve in $V_{\tau}$.
  Put $\bfu=\frac{d\xi(r)}{dr}(1)\in T_{\bfw}V_{\tau}$ and 
  let $\rho(z)=\sum_{j=1}^n |z_j|^2$.
  Note that
  \begin{eqnarray*}
   \bfu&=\frac{d\xi(r)}{dr}(1)
   =\frac{d\vphi_j(r,w_j)}{dr}|_{r=1}\,w_j\ne 0.
      \end{eqnarray*}
  Now we have
  \begin{eqnarray*}
   \begin{split}
    \frac{d\rho(\xi(r))}{dr}|_{r=1}&=\frac {d(\sum_{j=1}^n\eta(r,w_j)^2)}{dr}|_{r=1}\,\\
    &=2\sum_{j=1}^n\frac{d\vphi_j(r,w_j)}{dr}|_{r=1}\,\vphi(r,w_j)\,|w_j|^2>0
    \end{split}
\end{eqnarray*}
  This implies that $\bfu$ is not tangent to the sphere
  $S_{r}^{2n-1}$.
  \end{proof}

Fix a positive number $r$.
Choose a positive number $\eta_0>0$ so that 
the fibers $f_t\inv(\eta)$  and the sphere $S_r^{2n-1}$ intersect
 transversely for any $\eta,\,|\eta|\le \eta_0$.
Let 
\[
 \begin{split}
&\partial\mathcal
E(\eta_0,r):=\{(\bfz,t)\in \BC^n\times I\,|\,|f_t(\bfz)|=\eta_0,\,\|\bfz\|\le r\}\\
&\partial
E_t(\eta_0,r):=\{\bfz\in \BC^n\,|\,|f_t(\bfz)|=\eta_0,\,\|\bfz\|\le r\}.
\end{split}
\]
Note that $f_t:\partial
E_t(\eta_0,r)\to S_{\eta_0}^1$ is equivalent to the Milnor fibration of
$f_t$
by the second description (see \cite{OkaMix}).
 Using the Ehresmann's fibration theorem \cite{Wolf1}
to the projection:
\[
 \pi:S_r^{2n-1}\times I\to I,\quad
\pi':\partial\mathcal
E(\eta_0,r)\to I,
\]
 we obtain: % the main result:
 \begin{MainTheorem}{\rm(Isotopy Theorem)}\label{MainTheorem}
\begin{enumerate}
\item
  There exists an isotopy $h_t:S_r^{2n-1}\to S_r^{2n-1}$
  such that
  $h_0=\id$ 
and $ f_t(h_t(\bfz))=f_0(\bfz)$ for any $\bfz\in  S_r^{2n-1}$
 with $|f_0(\bfz)|\le \eta_0$ and 
  $h_t(K_0)=K_t$ for each $t,\,0\le t\le 1$ where $K_t=V_t\cap S_r^{2n-1}$.

 \item The Milnor fibrations of $f_{\bfa,\bfb}$ and $f_{\bfa}$ 
by the second description
\[
 f_0:\partial E_0(\eta_0,r)\to S_{\eta_0}^1,\,\,
f_t:\partial E_t(\eta_0,r)\to S_{\eta_0}^1
\] are
  $C^\infty$ equivalent.
\end{enumerate}
  \end{MainTheorem}
\begin{Remark}
As the first and the second description of Milnor fibrations
are equivalent (\cite{OkaMix}), the Milnor fibrations of $f_{\bfa,\bfb}$ and $f_{\bfa}$
 are
$C^\infty$ equivalent.
\end{Remark}
Applying the above  method to construct norm %and the argument
 preserving  diffeomorphisms
$h_t:\BC^n\setminus\{O\}\to \BC^n\setminus\{O\}$: $\|h_t(\bfz)\|=\|\bfz\|$, 
%$\arg(f_0(\bfz,\bar\bfz))=\arg(f_t(h_t(\bfz),\bar h_t(\bfz)))$
and
$h_t(V_{\bfa,\bfb})\subset V_t$, we obtain:
  \begin{Corollary}
   The pair $(\BC^n,V_t)$ is homeomorphic to the
   pair $(\BC^n,V_{\bfa,\bfb})$. This homeomorphism can be diffeomorphic 
outside
   of the origin. 
%This is a defeomorphism of the global Milnor fibrations.
   \end{Corollary}
\subsection{Mixed polar homogeneous projective hypersurfaces} Let us consider mixed
  homogeneous hypersurfaces case. Namely 
$a:=a_1=\dots=a_n$.
% and $b:=b_1=\dots=b_n$.
Thus $f_{\bfa,\bfb}=z_1^{a+b_1}\bar z_1^{b_1}+\cdots z_n^{a+b_n}\bar z_n^{b_n}$.
The hypersurface
$V_t:=\{f_t(\bfz,\bar\bfz)=0\}$ does not have $\BC^*$-action, as 
the polynomial $f_t$ is not radially homogeneous.
However $f_t$ is polar homogeneous. Thus $V_t$ has the canonical
$S^1$-action, defined by
$\la\circ \bfz=(\la z_1,\dots,\la z_n)$ for $\la\in S^1\subset \BC^*$.
Thus putting $K_t=V_t\cap S^{2n-1}$,
 the following diagram makes a good sense
\[
 \begin{matrix}
S^{2n-1}&\supset &K_t\\
\mapdown{\pi}&&\mapdown{\pi}\\
\BP^{n-1}&\supset &H_t
\end{matrix}
\]
where $H_t$ is the quotient space  $K_t/S^1$.
%Constructing the isotopy $h_t:S^{2n-1}\to S^{2n-1}$
%to be $S^1$-equivariant,
Using this family, we have the following result.
\begin{Corollary}
The projective hypersurfaces %$H_{\bfa,\bfb}
$H_{\bfa,\bfb}=\{[\bfz]\in \BP^{n-1}\,|\,f_{\bfa,\bfb}(\bfz,\bar\bfz)=0\}$
and $H_{\bfa}=\{[\bfz]\in \BP^{n-1}\,|\,f_{\bfa}(\bfz)=0\}$ are ambient 
isotopic. That is, there exists an isotopy 
$\bar h_t:\BP^{n-1}\to \BP^{n-1}$
%$\bar h_1: \BP^{n-1}\to \BP^{n-1}$
such that $\bar h_1(H_{\bfa,\bfb})=H_{\bfa}$.
\end{Corollary}
\section{Other polar weighted homogeneous polynomials}
In this section, we will generalize the previous results for simplicial
polar
weighted homogeneous polynomials which have isolated singularities at
the origin and whose  associated Laurent  polynomials are
simplicial weighted homogeneous polynomials
listed in \cite{OrlikWagreich:1971}.
%As an example, 
\subsection{Coefficients can be 1}
A mixed  polynomial $f(\bfz,\bar\bfz)=\sum_{i=1}^m
\,c_i\,\bfz^{\nu_i}\bar\bfz^{\mu_i}$ ($c_1,\dots, c_n\ne 0$)
is called simplicial if $m=n$ and the matrices $N\pm M$
are non-degenerate. Here the multi-integers are considered as column
vectors and the matrices $M,N$ are defined as
$N=(\nu_1,\dots, \nu_n)$ and $M=(\mu_1,\dots, \mu_n)$.
First we prove the following simple lemma.
\begin{Lemma}\label{canbe1}
Put $\tilde f(\bfw,\bar\bfw)=\sum_{i=1}^n \bfw^{\nu_i}\bar\bfw^{\mu_i}$.
Then there is a scaling linear mapping\,\,
$\vphi:\BC^n\to\BC^n$
defined by $\bfz\mapsto (w_1,\dots,w_n)=(z_1\alpha_1,\dots,
 z_n\alpha_n)$ such that  $\tilde f(\vphi(\bfz))=f(\bfz)$.
\end{Lemma}
\begin{proof}
Write $c_j=\exp(a_j+ib_j)$ and   $\alpha_j=\exp(\ga_j+i\eps_j)$ with
 $a_j,b_j,\ga_j,\eps_j\in \BR$ and
 $j=1,\dots, n$.
Then we consider the equality $c_j\bfz^{\nu_j}\bar\bfz^{\mu_j}=
\bfw^{\nu_j}\bar\bfw^{\mu_j}$ with $\bfw=\vphi(\bfz)$, $w_j=\al_j\,z_j$
 which reduces to:
\begin{eqnarray}\label{eq1}
  \al_1^{\nu_{j1}}\cdots \al_n^{\nu_{jn}}\bar\al_1^{\mu_{j1}}\cdots
\bar\al_n^{\mu_{jn}}=c_j,\, j=1,\dots, n.
\end{eqnarray}
The equality (\ref{eq1}) can be split in to the argument part and the
 absolute value  part as follows.
\[
 \begin{split}
&(\eps_1,\dots,\eps_n)(N-M)=(b_1,\dots, b_n)\\
&(\ga_1,\dots,\ga_n)(N+M)=(a_1,\dots, a_n)
\end{split}
\]
By the assumption, $\det(N\pm M)\ne 0$ and the scaling complex numbers
 $\al_1,\dots,\al_n$ are 
uniquely determined.
\end{proof}

\subsection{Other simplicial polar weighted polynomials and the
  generalization of the isotopy theorem}
We consider the following two simplicial polynomials with coefficient 1
where $g_1(\bfz),\,g_2(\bfz)$ are polynomials listed in
\cite{OrlikWagreich:1971}.
%An arbitrary weighted homogeneous simplicial polynomial with an isolated 
%singularity at the origin
\[
\begin{split}
&I:\begin{cases}
& f_I(\bfz,\bar\bfz)=z_1^{a_1+b_1}\bar
 z_1^{b_1}z_2+\cdots+z_{n-1}^{a_{n-1}+b_{n-1}}\bar
 z_{n-1}^{b_{n-1}}z_n+z_n^{a_n+b_n}\bar z_n^{b_n}\\
& g_I(\bfz)=z_1^{a_1}z_2+\cdots+ z_{n-1}^{a_{n-1}}z_n+z_n^{a_n}
\end{cases}\\
&
II:\begin{cases}
& f_{II}(\bfz,\bar\bfz)=z_1^{a_1+b_1}\bar
 z_1^{b_1}z_2+\cdots+z_{n-1}^{a_{n-1}+b_{n-1}}\bar
 z_{n-1}^{b_{n-1}}z_n+z_n^{a_n+b_n}\bar z_n^{b_n}z_1\\
& g_{II}(\bfz)=z_1^{a_1}z_2+\cdots+ z_{n-1}^{a_{n-1}}z_n+z_n^{a_n}z_1
\end{cases}
 \end{split}
\]
where $a_j\ge 1$ and $b_j\ge 0$ for each $j=1,\dots, n$.

Note that an arbitrary simplicial weighted homogeneous polynomial $g(\bfz)$ with
an isolated singularity at the origin is written as  joins
 of several simplicial
weighted homogeneous polynomials of either a Brieskorn type, or $g_I$ or
$g_{II}$.
To show the isotopy theorem for an arbitrary simplicial polar weighted homogeneous
polynomial,
it is enough to show the assertion for these three class of weighted
homogeneous polynomials. In the following section, we present the proofs
of the
 similar assertion
except for the transversality theorem for $f_{II}$ which is still open.
See Problem-Conjecture \ref{Conj} below.

\vspace{.3cm}
Consider the family of mixed hypersurfaces
$V_{\iota,t}=f_{\iota,t}\inv(0)\subset \BC^n$
where $\iota=I,II$ and 
$f_{\iota,t}(\bfz,\bar\bfz):=\,(1-t)f_\iota(\bfz,\bar\bfz)\,+\,t\,g_\iota(\bfz)$
for $0\le t\le 1$.
We investigate the similar assertions as in \S 2. First we have:
\begin{Lemma}The mixed polynomial function 
$f_{\iota,t}(\bfz,\bar\bfz):\BC^n\to \BC$
has %variety $V_{\iota,t}$ has 
a unique mixed
 singularity at
 the origin %and $V_{\iota,t}\setminus\{O\}$ is non-singular 
for any
 $t,\,0\le t\le 1$ and $\iota=I,II$. 
\end{Lemma}
\begin{proof}
The assertion is true for $t=0, 1$. Thus we assume that $0<t<1$.
Assume that $\bfw\in \BC^n\setminus \{O\}$ is a mixed singular point
of the function $f_{\iota,t}$.
Then we have by Proposition 1 of \cite{OkaPolar},
\begin{eqnarray}\label{a1}
 \overline{df_{\iota,t}(\bfw,\bar\bfw)}=\lambda \,\bar d\,
  f_{\iota,t}(\bfw,\bar\bfw),\quad\exists \la,\, |\la|=1.
\end{eqnarray}
Case 1. $\iota=I$. First assume that $w_n\ne 0$. 
Let $s=\min\{j\,|\, w_k\ne 0,\,k\ge j\}$.
Then by (\ref{a1}), we have 
\[
 \overline{\frac{\partial f_{\iota,t}}{\partial z_s}}(\bfw)=\la\,
 \frac{\partial f_{\iota,t}}{\partial \bar z_s}(\bfw).
\]
This gives the equality:
\begin{eqnarray*}
& {\bar w_n}^{a_n}\{(1-t)(a_n+b_n)|w_n|^{2b_n}+a_n\,t\}=\la\,
b_n w_n^{a_n} |w_n|^{2b_n}(1-t),
\,\text{if}\,\, s=n
\end{eqnarray*}
and if $s<n$, it gives:
\[
 {\bar w_s}^{a_s}\bar w_{s+1}\{(1-t)(a_s+b_s)|w_s|^{2b_s}+ a_s\,t\}=\la\,
b_s{w_s}^{a_s} w_{s+1}| w_s|^{2b_s}(1-t).\]
Denote  the left and the right side of the above equality
 by $L(s)$ and $R(s)$ respectively.
In the  both cases, we have a contradiction $|L(s)|> |R(s)|$ as
\[\begin{split}
& |L(s)|\ge\begin{cases}
	   (a_n+b_n)|w_n|^{a_n+2b_n}(1-t), %>b_n|w_n|^{a_n+2b_n}=|R(s)|,
\quad &s=n\\
   %   |L(s)|\ge      
  (a_s+b_s)|w_s|^{a_s+2b_s}|w_{s+1}|(1-t),\,\,&s<n
%>b_s       |w_s|^{a_s+2b_s}|w_{s+1}|(1-t)=|R(s)|.   
\end{cases}\\
&|R(s)|=\begin{cases}
b_n|w_n|^{a_n+2b_n}(1-t),\,\,&s=n\\
b_s       |w_s|^{a_s+2b_s}|w_{s+1}|(1-t),\,& s<n
\end{cases}.
\end{split}
\]
Nex we consider the case $w_n=0$. Let $\ell=\min\{j\,|\, w_k= 0,\,k\ge
 j\}$.
Then $\ell>1$. Consider the equality
\[
 \overline{\frac{\partial f_{\iota,t}}{\partial z_\ell}}{\bfw}=\la \frac{\partial
 f_{\iota,t}}{\partial \bar z_\ell}(\bfw).
\]
This gives an contradiction:
\[
\bar w_{\ell-1}^{a_{\ell-1}}\{|w_{\ell-1}|^{2b_{\ell-1}}(t-1)+t\}
\,=\, 0.
\]
Case 2. $\iota=II$.
Assume that $w_j=0$ for some $j$. Then we may assume that $w_n=0$ after
 the cyclic permutation of the index $k\mapsto k+n-j$ modulo $n$.
Then the proof is the exact same as in Case 1 with $w_n=0$.

Assume that $\bfw\in \BC^{*n}$.
Then for each $j$, we have
\begin{multline}\label{eq2}
 {\bar w_s}^{a_s}\bar w_{s+1}\{(1-t)(a_s+b_s)|w_s|^{2b_s}+ a_s\,t\}+\\
(1-t)\bar w_{s-1}^{a_{s-1}}\bar w_{s}|w_{s-1}|^{2b_{s-1}}
=\la\,
b_s{w_s}^{a_s} w_{s+1}| w_s|^{2b_s}(1-t).
\end{multline}
Here the numbering is understood modulo $n$,  so $w_j=w_{j+n}$ etc.
Consider an index  $m$  so that $|w_m|^{a_m+2b_m}
 |w_{j+1}|\ge|w_j|^{a_j+2b_j} |w_{j+1}|$ for  any $j$.
Let $L(m)$ and $R(m)$ be the left and right side quantities of
 (5)
for $s=m$.
Then 
\[\begin{split}
  & |L(m)|>\\
&(1-t)(a_m+b_m)|w_m|^{a_m+2b_m}|w_{m+1}|-
(1-t)|w_{m-1}|^{a_{m-1}+2b_{m-1}}|w_m|\\
&\ge (1-t)(a_m+b_m-1)|w_m|^{a_m+2b_m}|w_{m+1}|\\
&\ge (1-t)b_m
 |w_m|^{a_m+2b_m}|w_{m+1}|
=|R_m|
\end{split}
 \]
which is an contradiction to (\ref{eq2}). This completes the proof.
\end{proof}
\newtheorem{MainTheoremBis}[Theorem]{Main Theorem-bis}\label{Conj}%
The next key Lemma is:
\begin{Lemma}\label{a3} For any $0\le t\le 1$ and $r>0$,
the sphere $S_r^{2n-1}$ intersects transversely with the mixed
 hypersurface
$V_{I,t}$.
\end{Lemma}
\begin{proof}
The proof of Lemma \ref{a3} is  more complicated as that of 
Lemma \ref{key2}.  We assume that $\iota=I$.
Take a point $\bfw\in V_{\iota,t}\setminus\{O\}$. Put $\rho=\|\bfw\|$.
For the proof, it suffices to show that
$T_{\bfw}V_{\iota,t}\not \subset T_{\bfw}S_{\rho}^{2n-1}$.
Consider first the sets 
%$I_0=\{i\,|\,w_i=0 \}$
%
%Case 1. Assume that $\iota=I$.
% and put
\[
I_0=\{i\,|\,w_i=0 \},\quad  J=\{j\,|\, w_j^{a_j}w_{j+1}^{\eps_{j,n}}\ne 0\}
%&\iota=I\\
%\{j\,|\, w_j^{a_j}w_{j+1}\ne 0\},\,\,&\iota=II
%\end{cases}
\]
where $\eps_{j,n}=1$ for $1\le j<n$ and $0$ for $j=n$.
%Note that $j\in J$ if and only if $w_j,w_{j+1}\ne0$ where $w_{n+1}:=1$ 
%by definition.

 Assume that $J=\emptyset$. Then
we put 
$\bfw_j(s)=\bfw_j$ for $j\in I$ and $\bfw_j(s)=s\bfw_j$ for $s\notin I_0$.
Then it is easy to see that 
$f_{\iota,t}(\bfw(s),\bar \bfw(s))\equiv 0$ for $s> 0$ and 
$\|\bfw(s)\|$ is obviously   monotone increasing in $s$.
Thus the tangent vector $\bfv:=\frac{d\bfw(s)}{ds}|_{s=1}$ is contained
 in
$T_{\bfw}V_{\iota,t}\setminus T_{\bfw}S_{r}^{2n-1}$.

Assume that $J\ne 0$.
A subset $K\subset J$ is called {\em connected} if 
$i,j\in K,\,i<j$ implies   $k\in K$ for any $k,\,i\le k\le j$.
%We may assume that $1\in I$  if $I\ne \emptyset$ and $\iota=II$, $1\in
 %I$.

Case 1. %We assume either $\iota=I$ or $\iota=II$ and $1\in I$.
Let $J=J_1\cup \cdots \cup J_\ell$ be the decomposition into the
 connected components of $J$. We may assume that
$J_i=\{k\in \BN\,|\,\nu_j\le k\le \mu_j\}$  with some 
$\nu_i\le \mu_i\in J_i$
and 
\[
 \mu_i+1<\nu_{i+1},\quad i=1,\dots, \ell-1.
\]

We consider the following system of equations for positive parameters
$s_1,\dots, s_n$ for a given $r>0$. Put 
$z_j(s_j)=w_js_j$ for $j=1,\dots,n$.
We are going to show the existence of functions $s_j=s_j(r),\,j\in J_i$
 so that 
\[
 f_t(z_1(s_1),\dots, z_n(s_n),\bar z_1(s_1),\dots, \bar z_n(s_n))\equiv
 0,
\quad \frac{d\|\bfz(\bfs)\|}{dr}>0.
\]
First we put $ s_j(r)\equiv 1$ for $ j\notin J$.
We fix $i$ and we will show that there exists a differentiable 
solution 
$(s_{\nu_i}(r),\dots, s_{\mu_i}(r))$ of the equations:
\[\begin{split}
 {\rm E_j}:\, 
z_j(s_j)^{a_j}z_{j+1}(s_{j+1})&\{|z_j(s_j)|^{2b_j}(1-t)+t\}\\
 \qquad &=r\,
   w_j^{a_j}w_{j+1}\{|w_j|^{2b_j}+t\},
\,\, \text{for}\,\,\,\nu_j\le j\le \mu_j<n
\end{split}
\]
For the induction purpose, we rewrite this equation  as follows.
\[\begin{split}
 {\rm E_j'}:\, 
w_j^{a_j}s_j^{a_j}w_{j+1}&\{|z_j|^{2b_j}s_j^{2b_j}(1-t)+t\}\\
 \qquad &=\frac{r}{s_{j+1}}\,
   w_j^{a_j}w_{j+1}\{|w_j|^{2b_j}+t\},
\,\, \text{for}\,\,\,\nu_j\le j\le \mu_j<n
\end{split}
\]
For the case $j=\mu_i=n$, this $E_n$ takes the form:
\[
  {\rm E_n}:\, z_n^{a_n}s_n^{a_n}\{|z_n|^{2b_n}s_n^{2b_n}(1-t)+t\}
=
   r\, w_n^{a_n}\{|w_n|^{2b_n}+t\}.
\]
%We fix $i$ and we will show that there exists a differentiable 
%solution \nl
%$(s_{\nu_i}(r),\dots, s_{\mu_i}(r))$.

%and ${\bf e}:=(1,\dots, 1)$.
We solve the equality using the same argument as in the proof of Lemma
 3 from $j=\mu_i$ to $j=\nu_i$ downward.
For this purpose, we use the equality $E_j'$ by replacing $r/s_{j+1}$
 by $r_{j}$:
\[\begin{split}
 E_j'':\,\,w_j^{a_j}s_j^{a_j}w_{j+1}&\{|w_j|^{2b_j}s_j^{2b_j}(1-t)+t\}\\
&\qquad =r_j\,
   w_j^{a_j}w_{j+1}^{\eps_{j+1,n}}\{|w_j|^{2b_j}+t\},\, %r_j>0.
\,\nu_j\le j\le \mu_j.
\end{split}
\]
First we solve $E_{j}''$. Put $\phi_j(s_j)$ be the left side quantity of $E_j''$.
By the monotone property of the real valued function
 $|\phi_{j}(s_{j})|$
and by the property 
$$\arg(z_{j}^{a_j}(s_j)z_{j+1}(s_{j+1}))=\arg(w_{j}^{a_j}w_{j+1})=\text{constant}$$
%and $\phi_{j}(\bfe)=|w_{j}|^{a_{j}}\{(1-t)|w_{j}|^{2b_{j}}+t\}$,
we can solve $s_j$ as a positive valued  differentiable function of
 $r_j$, say $s_j=\psi_j(r_j)$,
so that 
%$\psi_j(1)= 1$ and 
\begin{eqnarray}\label{induction}
\psi_j(1)= 1,\,\, \frac{d\psi_{j}}{dr_j}(r_j)>0,\,\,
\psi_j(r_j)^{a_j}\le r_j.
\end{eqnarray}
Now we are ready to solve $E_j',\,j\in J_i$ inductively from $j=\mu_i$.
First we put $r_{\mu_i}=r$  as $s_{\mu_i+1}=1$ and 
$s_{\mu_j}(r)=\psi_{\mu_i}(r)$.
Then inductively we define $r_j(r),\,s_j(r)$ by
\[\begin{split}
% &r_{\mu_i-1}(r):=r/s_{\mu_i}(r),
%\, s_{\mu_i-1}(r):=\psi_{\mu_i-1}(r_{\mu_i-1}(r)),\\
&r_j(r):=r/s_{j+1}(r),\,
s_j(r):=\psi_j(r_j(r)),\,\,
j=\mu_i-1,\dots,\nu_i.
\end{split}
\]
Note that $s_j(r),\,j=1,\dots, n$ are real valued functions and
 differentiable in $r$.
%It is easy to observe that  $E_j'$ implies  $E_j$ by substituting
%$r_j=r/s_{j+1}$.
By the inequality (\ref{induction}),
\[
 r_{\mu_i-1}(r)=r/s_{\mu_i}(r)\ge r^{1-1/a_{\mu_i}}\ge 1,
\, s_{\mu_i-1}(r)\ge 1
\quad \text{for}\,\,r\ge 1.
\]
We show by induction
that
\[
 r_j(r)\ge 1,\,\,s_j(r)\ge 1\quad \text{for}\,\, r\ge 1,\,\nu_i\le
 j\le \mu_i.
\]
This is true for $j=\mu_i$ and $\mu_i-1$ as we have seen above. For $j\le \mu_i-2$,
\[\begin{split}
 r_j(r)^{a_{j+1}}&=\frac{r^{a_{j+1}}}{s_{j+1}(r)^{a_{j+1}}}=
\frac{r^{a_{j+1}}}{\psi_{j+1}(r_{j+1}(r))^{a_{j+1}}}\underset{(\ref{induction})}\ge
\frac{r^{a_{j+1}}}{r_{j+1}(r)}\\
&\ge {r^{a_{j+1}-1}}\cdot{s_{j+2}(r)}\ge  1\quad
   \text{for}\,\nu_i\le j\le\mu_i-2,\,r\ge 1\\
\end{split}
\]
This implies $r_j(r)\ge 1$ 
and $s_j(r)\ge 1$ for $r\ge 1$.

Now we observe that %Thus for $r\ge 1$,
\[
 \|(s_1(r)w_1,\dots, s_n(r) w_n)\|\ge 
\|\bfw\|,\quad \text{for}\,\,r\ge 1
\]
and 
\[
\left \|\frac{d(s_1(r)w_1,\dots, s_n(r) w_n)}{dr}\right\|>0
\]
as we have 
\[
 \frac{d(|s_{\mu_i}(r)w_{\mu_i}|}{dr}>1.
\]
Now by taking the summation of $E_j$ for $j\in J_i$ and
 $i=1,\dots,\ell$,
we get the equality:
\[
 f_{\iota,t}(\bfz(r))=r\,f_{\iota,t}(\bfw)\equiv 0\quad
\text{where}\,\,
\bfz(r)=(s_1(r)w_1,\dots, s_n(r) w_n).
\]
Let $\bfv=\frac{d\bfz(r)}{dr}|_{r=1}$.
Then $\bfv\in T_{\bfw} V_{\iota,t}\setminus  T_{\bfw}S_{\rho}^{2n-1}$
by the above observation.
This completes the proof.
\end{proof}
\begin{Remark}
The above proof does not work if $J=\{1,\dots,n\}$ and $\iota=II$,
as we do not have a starting point of the induction.
\end{Remark}

\noindent
\newtheorem{ProbConj}[Theorem]{Problem-Conjecture}
\begin{ProbConj}
Is the assertion  true for $V_{II,t}$?
 \end{ProbConj}

Now we assert that 

\begin{MainTheoremBis}\label{MainTheoremBis}
Fix a positive $r$.
 Choose a sufficiently small positive $\eta_0$
so that $f_{I,t}\inv(\eta)$ and $S_r^{2n-1}$ intersect transversely for any 
$\eta,\,|\eta|\le \eta_0$.
\begin{enumerate}
\item There exists an isotopy
$h_{t}:(S_r^{2n-1},K_{0,r})\to (S_r^{2n-1},K_{t,r})$
such that 
$f_{I,t}(h_t(\bfz))=f_{I,0}(\bfz)$ for any $\bfz$ with $|f_{I,0}(\bfz)|\le \eta_0$.
\item 
The Milnor fibrations
\[
 f_{I,0}:\partial E_0(\eta_0,r)\to S_{\eta_0}^1,\,\,
f_{I,t}:\partial E_t(\eta_0,r)\to S_{\eta_0}^1
\] are
$C^\infty$ equivalent 
where 
\[
 \partial E_t(\eta_0,r)=\{\bfz\,|\,
      |f_{I,t}(\bfz)|=\eta_0,\,\|\bfz\|\le r\}.
\]
\end{enumerate}
\end{MainTheoremBis}
The assertion is also true for $f_{II}$ if the above transversality conjecture is true.

\begin{Corollary}
The Milnor fibrations of $f_I(\bfz,\bar\bfz)$ and $g_I(\bfz)$ are
 $C^\infty$ equivalent.
\end{Corollary}

%\bibliographystyle{abbrv}
%\bibliography{/home/oka/paper/okabib}

\begin{thebibliography}{1}

\bibitem{OkaPolar}
M.~Oka.
\newblock Topology of polar weighted homogeneous hypersurafces.
\newblock {\em Kodai Math. J.}, 31(2):163--182, 2008.

\bibitem{OkaMix}
M.~Oka.
\newblock Non-degenerate mixed functions.
\newblock{arXive:  math 0909.1904},
\newblock {\em to appear Kodai Math. J.}

\bibitem{OrlikWagreich:1971}
P.~Orlik and P.~Wagreich.
\newblock Isolated singularities of algebraic surfaces with {C}$\sp{\ast}$\
  action.
\newblock {\em Ann. of Math. (2)}, 93:205--228, 1971.

\bibitem{R-S-V}
M.~A.~S. Ruas, J.~Seade, and A.~Verjovsky.
\newblock On real singularities with a {M}ilnor fibration.
\newblock In {\em Trends in singularities}, Trends Math., pages 191--213.
  Birkh\"auser, Basel, 2002.

\bibitem{Wolf1}
J.~A. Wolf.
\newblock Differentiable fibre spaces and mappings compatible with {R}iemannian
  metrics.
\newblock {\em Michigan Math. J.}, 11:65--70, 1964.

\end{thebibliography}
\def\cprime{$'$} \def\cprime{$'$} \def\cprime{$'$} \def\cprime{$'$}
  \def\cprime{$'$} \def\cprime{$'$} \def\cprime{$'$} \def\cprime{$'$}

\end{document}